\documentstyle{article}

\catcode`\@=11
\font\tenmsx=msxm10
\font\sevenmsx=msxm7
\font\fivemsx=msxm5
\font\tenmsy=msym10
\font\sevenmsy=msym7
\font\fivemsy=msym5
\newfam\msxfam
\newfam\msyfam
\textfont\msxfam=\tenmsx  \scriptfont\msxfam=\sevenmsx
 \scriptscriptfont\msxfam=\fivemsx
\textfont\msyfam=\tenmsy  \scriptfont\msyfam=\sevenmsy
 \scriptscriptfont\msyfam=\fivemsy
    
\def\hexnumber@#1{\ifnum#1<10 \number#1\else
  \ifnum#1=10 A\else\ifnum#1=11 B\else\ifnum#1=12 C\else
  \ifnum#1=13 D\else\ifnum#1=14 E\else\ifnum#1=15 F\fi\fi\fi\fi\fi\fi\fi}
\def\msx@{\hexnumber@\msxfam}
\def\msy@{\hexnumber@\msyfam}
      
\mathchardef\upharpoonright="3\msx@16     
\let\restriction=\upharpoonright
\mathchardef\Vdash="3\msx@0D
\mathchardef\Vvdash="3\msx@0E
\mathchardef\vDash="3\msx@0F
\mathchardef\nleq="3\msy@02
\mathchardef\ngeq="3\msy@03
\mathchardef\nless="3\msy@04
\mathchardef\ngtr="3\msy@05
\mathchardef\nprec="3\msy@06 
\mathchardef\nsucc="3\msy@07   
\mathchardef\lneqq="3\msy@08
\mathchardef\gneqq="3\msy@09  
\mathchardef\nleqslant="3\msy@0A
\mathchardef\ngeqslant="3\msy@0B
\mathchardef\lneq="3\msy@0C
\mathchardef\gneq="3\msy@0D
\mathchardef\npreceq="3\msy@0E   
\mathchardef\nsucceq="3\msy@0F 
\def\Bbb{\ifmmode\let\next\Bbb@\else
 \def\next{\errmessage{Use \string\Bbb\space only in math mode}}\fi\next}
 \def\Bbb@#1{{\Bbb@@{#1}}}
 \def\Bbb@@#1{\fam\msyfam#1}
\catcode`\@=\active

\newtheorem{theorem}{Theorem}[section]
\newtheorem{defn}[theorem]{Definition}
\newtheorem{prop}[theorem]{Proposition}

\newtheorem{cor}[theorem]{Corollary}
\newtheorem{claim}[theorem]{Claim}
\newtheorem{eple}[theorem]{Example}

\newfont{\gothic}{eufm10}

\newcommand{\implies}{\rightarrow}
\newcommand{\sub}{\subseteq}

\newcommand{\sm}{\setminus}
\newcommand{\proof}{{\bf Proof:}\ }
\newcommand{\qed}{\hfill $\Box$}
\newcommand{\noi}{\noindent}
\newcommand{\rest}{\restriction}
\newcommand{\ra}{\rangle}
\newcommand{\la}{\langle}

\newcommand{\CA}{{\cal A}}
\newcommand{\CB}{{\cal B}}

\newcommand{\CF}{{\cal F}}
\newcommand{\CG}{{\cal G}}
\newcommand{\CH}{{\cal H}}

\newcommand{\CM}{{\cal M}}
\newcommand{\CP}{{\cal P}}
\newcommand{\CS}{{\cal S}}
\newcommand{\CT}{{\cal T}}
\newcommand{\CU}{{\cal U}}
\newcommand{\CV}{{\cal V}}

\newcommand{\dCF}{{\cal F}\!\!\downarrow}

\newcommand{\gb}{{\gothic b}}
\newcommand{\gd}{{\gothic d}}
\newcommand{\gc}{{\gothic c}}
\newcommand{\gu}{{\gothic u}}
\newcommand{\gothg}{{\gothic g}}
\newcommand{\dgb}{{\gothic b}$\downarrow$}

\newfont{\speck}{msam5}
\newcommand{\oo}{\omega^{\mbox{{\tiny $\uparrow$}}\omega}}
\newcommand{\so}{[\omega]^\omega}
\newcommand{\fa}{\forall}
\newcommand{\ex}{\exists}

\newcommand{\exi}{\exists^\infty}
\newcommand{\bd}{\begin{description}}
\newcommand{\ed}{\end{description}}

\title{Bounding and dominating number \\
        of families of functions on N}
     
\date{}
\author{Claude Laflamme\thanks{This research was partially supported by 
         NSERC of Canada.  \newline
        1980 {\em AMS Subj. Class. (1985 Revision)} Primary 03E35; Secondary 04A20.
        } \\
        Department of Mathematics and Statistics \\
        University of Calgary \\
        Calgary, Alberta  \\
        Canada T2N 1N4}

\begin{document}
     
\maketitle

\begin{abstract}

We pursue the study of families of functions on the natural numbers,
with emphasis here on the bounded families.  The situation being more
complicated than the unbounded case, we attack the problem by
classifying the families according to their bounding and dominating
numbers, the traditional scheme for gaps.  Many open questions remain. 

\end{abstract}

\section{Introduction}

Over the years, the notion of gaps of functions (or sets) of natural
numbers has played an important role in the application of Set Theory to
different branches of mathematics, see for example \cite{Scheepers} for
a survey.  It should thus not come as a surprise that families of
functions in general might have an impact.  It was in fact shown
recently in \cite{Blass1} that the structure of directed unbounded
families of functions had an influence on several problems and in
\cite{JustLaf} an application of non-directed unbounded families was
made.  In the papers \cite{Blass1, BlassLaf, Laf1, Laf2}, a sufficiently
precise description of these families was given to address these
questions. 

The next step is to consider bounded families and provide a similar
description.  Such families are no more than generalizations of the
classical notion of (linearly ordered) gaps and the situation appears
quite complex; indeed, not only all partial orders of size $\leq
\aleph_1$ embedd as bounded families of functions but unbounded families
themselves reflect as bounded ones.  We have tried in this
paper to classify the bounded families according to their bounding and
dominating numbers, a criterion much weaker than cofinal equivalence for
example, but which seemed a good starting point but in fact many open
questions remain. 

The ultimate results we hope to achieve are to describe the families
that can be built from ZFC alone with enough details to be useful for
applications; forcing is used to verify that ZFC has been exhausted, i.e.
that no other families can be built from ZFC.

We express our warm thanks to Stevo Todorcevic for discussions on the topic.

\section{Notation and Preliminaries}

We write $\omega$ for the set of natural numbers, $\omega^\omega$ for
the set of all functions on $\omega$ and $\oo$ for the set of
non-decreasing (monotone) functions.  We often use $f \leq g$ to
abbreviate $(\forall n)f(n) \leq g(n)$, and $f \leq^{*} g$ for
$(\forall^{\infty}n)f(n) \leq g(n)$ and similarly for $<$ and $<^*$;
here ``$\forall^{\infty}n$" means ``for all but finitely many $n$" and
similarly ``$\exists^{\infty}n$" means ``there exists infinitely many
$n$".  Also important is the ordering $f \prec g$ defined by $\lim_n
g(n) - f(n) = + \infty$.  As $\langle \omega^\omega, \prec \rangle$
embedds in $\langle \oo, \prec \rangle$, we shall be interested
essentially in the latter structure and we shall assume for the
remaining of this paper that we deal with monotone functions only. 

We use $\cal F,G,H$ to denote families of
functions. Further, we shall assume that the families considered are closed 
downward under $\leq^*$; this simplifies greatly the discussion without interfering
with the results. In particular the bounding and dominating numbers that we 
define are the same for a given  family or its downward closure, and
further it suffices to present the generators to describe a family.

Given two families of  functions $\CH \sub \CF$, we say that
$\CH$ is unbounded in $\CF$ if:

\[(\fa f \in \CF)(\ex h \in \CH) \;  h \nless^* f,\]

\noi and $\CH$ is said to dominate $\CF$ if

\[(\fa f \in \CF)(\ex h \in \CH)\;  f \leq^* h.\]

\noi We shall in practice prove slightly more, for example that $f \prec g$
instead of only $f \leq^* g$ or that $\limsup_n f(n) - g(n) = + \infty$ 
instead of just $f \nless^* g$; the reason is that often for applications
two functions are identified if they differ only by a fixed natural
number.  Indeed we shall work mostly with the $\prec$ ordering.

\noi A family dominating $\omega^\omega$ is usually called a dominating
family, and one unbounded in $\omega^\omega$ is simply called an unbounded 
family.

We define the bounding number \gb$(\CF)$ of a family $\CF$ as

\[ \mbox{\gb}(\CF) = \min \{ |\CH|: \CH \sub \CF \mbox{ is unbounded in } \CF \} \]

\noi and the dominating number \gd$(\CF)$ as

\[ \mbox{\gd}(\CF) = \min \{|H|: \CH \sub \CF \mbox{ dominates } \CF \}. \]

\noi \gb=\gb$(\oo)$ is the usual bounding number and \gd=\gd$(\oo)$ the dominating
number.

 The infinite subsets of $\omega$ are denoted by $[\omega]^{\omega}$,
the standard ordering is $A \sub^* B$ if $A \sm B$ is finite.  We shall
be interested in almost disjoint families, that is families of infinite
sets with pairwise finite intersections.  By fixing a bijection from the
rationals and $\omega$ and considering for each irrational number a
sequence of rationals converging to it, we see that there is an almost
disjoint family of subsets of $\omega$ of size \gc, the continuum. 
Typical functions that will interest us are of the form $next(-,X)$ for
some $X\in [\omega]^\omega$, where

\[ next(n,X)= \mbox{ the smallest element of $X$ greater
than or equal to $n$}, \]

\noi and similarly for the function $last(-,X)$. 
     
An ultrafilter is a proper family of infinite sets closed under finite
intersections, supersets, and maximal with respect to those properties;
in particular it must contain $X$ or $\omega \setminus X$ for any $X
\subseteq \omega$, and must be nonprincipal. 
 We use $\CU,\CV$ to denote ultrafilters.  We write $\chi(\CU)$ for the
minimal cardinality of a collection generating the ultrafilter $\CU$,
and {\gothic u} for the least cardinality of a family of sets generating
any ultrafilter.  A $P_\kappa$-point is an ultrafilter $\CU$ with the
property that any $\kappa$ decreasing sequence from $\CU$ has a lower
bound {\it in} $\CU$.

\section{Unbounded Families}

We shall consider in this section three sorts of unbounded (closed
downward) families (of monotone functions). 

\begin{defn}\label{unb-fam} ~ 
     
\begin{enumerate}
  \item The $\cal D$-class (the dominating class):
 $ {\cal F} \in {\cal D} \mbox{ iff } {\cal F} \mbox{ is dominating}.$
     
  \item The $\cal S$-class (the superperfect class):
    $ {\cal F} \in {\cal S} \mbox{ iff }$
    \begin{description}
       \item[i)] $(\exists h)(\forall f \in {\cal F})
                  (\exists^{\infty}_{n})[f(n) \leq h(n)] $
       \item[ii)] $(\exists g)(\forall f)
                  [(\exists^{\infty}_{n}) f(n) \leq g(n)
                         \rightarrow f \in {\cal F}] $.
    \end{description}
  \item The $\cal U$-class (for an ultrafilter $\cal U $):
    $ {\cal F} \in {\cal U} \mbox{ iff }$
    \begin{description}
       \item[i)] $ (\exists h)(\forall f \in {\cal F})
                [ \{n: f(n) \leq h(n)\} \in {\cal U}]$.
       \item[ii)] $ (\exists g)(\forall f)
                 [ \{n: f(n) \leq g(n)\} \in {\cal U}
                           \rightarrow f \in {\cal F}] $.
     \end{description}
 \end{enumerate}

\end{defn}

These three classes are easily seen to be distinct and we have shown in
\cite{Laf1,Laf2} the relative consistency of any unbounded family of
functions falling into one one these classes; that is in ZFC alone, no
other unbounded family of functions can be obtained and even a single
ultrafilter of your choice may be used for all members of the
3$^{\mbox{{\scriptsize rd}}}$ class.  This fulfills our original
motivation for unbounded families, in other words these descriptions are
sufficiently detailed to provide answers to many general mathematical
problems (see \cite{Blass1, BlassLaf, JustLaf, Laf2}).

The reason to pursue their studies here is their influence on bounded
families as we will see in the next section.  So we now discuss the
bounding and dominating number of families in these three classes.  If
$\CF$ is dominating, then \gb$(\CF)=$\gb, the usual bounding number and
\gd$(\CF)=$\gd, the usual dominating number.  Although it is possible to
make structural distinctions between dominating families, applications
have not made them yet necessary to analyze.  If we demand that our
families be closed downward under $\leq^*$, then there is only one
domimnating family, namely $\oo$. 

\vspace{.2in}

We now turn to the $\CS$-class. In \cite{Kechris}, Kechris showed that any 
unbounded Borel family must contain a superperfect tree, and we showed in 
\cite{Laf1} that any non-dominating family containing a superperfect tree 
belongs to the $\CS$-class.

\begin{prop}\label{S-class}
If $\CF$ is in the $\CS$-class, then \gb($\CF$)$=1$ or $2$ and \gd$=$\gc.
Further these two values of \gb \/ are attainable.
\end{prop}

\proof

\noi {\bf \ref{S-class}.1:} We first show that \gb($\CF$)$\leq 2$ if $\CF$
belongs to the $\CS$-class. 

The point is that a bounding number of at least 3 means that the family
is directed; it thus suffices to show that a directed family sastisfying
2ii) in the $\CS$-class is dominating, a contradiction. 

So let $\CF$ belong to the $\CS$-class and witnessed by $g$ and $h$ as
in definition \ref{unb-fam}.  Fixing any $p \in \oo$, define a sequence
of integers by $\pi_0=0$ and more generally such that $g(\pi_{n+1}) >
p(\pi_n)$.  If we now let $X_i = \{\pi_{2n+i}: n \in \omega \}$ for
$i=0,1$ and define $f_i(n)=g(next(n,X_i)) \in \CF$, then
$\max\{f_0(n),f_1(n)\} \geq p(n)$ for each $n$.  Since $p$ was
arbitrary, we see that $\CF$ must be dominating if directed, i.e.  if
\gb($\CF$)$\geq 3$. 

\noi {\bf \ref{S-class}.2:} We now show that \gd$(\CF)=$\gc \/ if $\CF$ belongs
to the $\CS$-class. 

Again, let $\CF$ belong to the $\CS$-class and witnessed by $g$ and $h$. 
Choose also an increasing sequence of integers $X=\{\pi_n : n \in \omega
\}$ such that $g(\pi_{n+1}) > h(\pi_n)$ for each $n$.  Now for any
infinite $Y \sub X$, consider the function $f_Y(n)=g(next(n,Y))$; this
function must belong to $\CF$ as it is equal to $g$ infinitely often. 
Observe however that for $p \in \CF$ such that $f_Y \leq^* p$,
 
\[S_p= \{n : p(n) \leq h(n) \} \sub^* S(Y)= \bigcup\{(\pi_{n-1},\pi_n]: \pi_n \in Y \} .\]

\noi Since moreover $Y \cap Y' =^* \emptyset \implies S(Y) \cap S(Y') =^*
\emptyset$, choosing an almost disjoint family of infinite subsets of
$X$ of size \gc \/ shows that \gc \/ functions are necessary to dominate
$\CF$. 

\noi {\bf \ref{S-class}.3:} An example of $\CF$ in the $\CS$-class with \gb($\CF$)=1 (and
\gd($\CF$)=\gc). 

Fix any unbounded function $g \in \oo$ and let
$\CF= \{ f \in \oo : (\exi n) f(n) \leq g(n) \} = \{g(next(-,X)): X \in \so \}$.
Since $g$ itself is unbounded in $\CF$ (and belongs to $\CF$), we have  \gb($\CF$)=1.

\noi {\bf \ref{S-class}.4:} An example of $\CF$ in the $\CS$-class with \gb($\CF$)=2 (and
\gd($\CF$)=\gc). 

Consider the identity function $id(n)=n$ and for $X \in \so$ let 
   \[ h_X(n)= next(n,X) + |X^c \cap n| \]
 and finally put $\CF = \{h_X: X \in \so \}$.

The fact that $\CF$ belongs to the $\CS$-class is witnessed by the
functions $g=id$ and $h(n)=2n$.  Then \gb($\CF$)$\leq 2$ and it thus
suffices to show that no single member of $\CF$ is unbounded.  But given
$h_X \in \CF$, choose an infinite $Y \sub X$ such that $X \setminus Y$
is also infinite.  Then for each $N$ and $n$ large enough

$\begin{array}{ll}
h_Y(n) &=next(n,Y) + |Y^c \cap n| \\
       &\geq next(n,X) + |Y^c \cap n| \\
       &\geq next(n,X) + |X^c \cap n| + N \\
       &= h_X(n) + N.
\end{array}$

\noi Thus $h_X \prec h_Y$ and hence \gb($\CF$)=2.  The proof
of \ref{S-class}.2 will actually give you two specific functions
unbounded in $\CF$.  This completes the proof of Proposition \ref{S-class}.  \qed

\medskip

We now turn our attention to members of the $\CU$-class.

\begin{prop} \label{U-class}
For any ultrafilter $\CU$,  $\lambda \in \{1,2,\omega\}$
and $\chi(\CU) \leq \kappa \leq $\gc, there is a family $\CF$ in the
$\CU$-class such that \gb($\CF$)$=\lambda$ and \gd($\CF$)$=\kappa$. 
\end{prop}

\proof First choose two functions $g, h \in \oo$ and an increasing
sequence of integers $\la   \pi_n : n \in \omega\ra  $ such that:

\begin{enumerate}
\item $g \prec h$
\item $\lim_n h(n+1)-[h(n)+n]=+\infty$
\item $(\fa n)\; h(\pi_n) \geq g( \pi_{n+2})$
\item $(\fa n)\; \pi_{n+1} - \pi_n \geq 2^n.$
\end{enumerate}

\noi Also fix an almost disjoint family $\CA= \la   A_\alpha: \alpha < $\gc$\ra  $
such that 
\[|A_\alpha \cap [\pi_n, \pi_{n+1})|=1 \mbox{ for each $n$ and $\alpha$.} \]

\noi Without loss of generality we may assume that $E= \bigcup_n
[\pi_{2n},\pi_{2n+1}) \in \CU$ as the other case is analogous.  For $X
\in \CU$.  let $P(X)= \cup\{[\pi_{2n-1},\pi_{2n}): X \cap
[\pi_{2n},\pi_{2n+1}) \neq \emptyset$. 

\vspace{.2in}

With these preliminaries we are ready to build the desired families.

\noi {\bf \ref{U-class}.1:} We build an $\CF$ in the $\CU$-class such that \gb($\CF$)=1
and \gd($\CF$)$=\kappa$ where
           $\chi(\CU) \leq \kappa \leq $\gc.

\noi For $X \in \CU$ such that $X \sub E$ and any $\alpha < \kappa$, let

\[ f^X_\alpha(n)= h(next(n,[A^c_\alpha \cap P(X)] \cup X)) \]

\noi and put $\CF = \{f^X_\alpha : X \sub E, X \in \CU, \alpha < \kappa \}$. 
Observe that $g(next(-,X)) \leq f^X_\alpha$ for any $\alpha$ and that if
$X, \alpha$ are given, then $f^X_\alpha(n)=h(x)$ for any $x \in X$ and
therefore $g,h$ witness that $\CF$ belongs to the $\CU$-class. 

As $h$ is unbounded in $\CF$, we conclude readily that \gb($\CF$)=1.  We
must now show that the dominating number is $\kappa$.  If $\CB$ is a base
for the ultrafiler $\CU$, then

\[\CH= \{f^X_\alpha : X \in \CB, \alpha < \kappa \} \]

\noi clearly dominates the family $\CF$ and therefore \gd($\CF$)$\leq
\chi(\CU) \cdot \kappa =\kappa$.  On the other hand, fix a family $\CH
\sub \CF$ of size less than $\kappa$, and fix some ordinal $\beta
\in\kappa$ not mentionned in any indexing of the functions from $\CH$. 
But if $X \in \CU, X \sub E$ and $\alpha \neq \beta$, then $f^X_\alpha$
does not dominate $f^E_\beta$.  Indeed the set $A^c_\alpha \cap P(X)
\setminus A^c_\beta$ is infinite as otherwise $P(X) \cap A_\beta \sub^*
A_\alpha$, and since $P(X)$ is an infinite union of intervals of the
form $[\pi_n,\pi_{n+1})$, $P(X) \cap A_\beta $ is infinite and thus
$A_\beta \cap A_\alpha$ is infinite contradicting that $\CA$ is an
almost disjoint family. 

But now for any $x \in A^c_\alpha \cap P(X) \setminus A^c_\beta$,
$f^X_\alpha (x)=h(x)$ and $f^E_\beta (x) \geq h(x+1)$ and as $\lim_n
h(n+1)-h(n) = +\infty$ we get $\limsup_n f^E_\beta (n) - f^X_\alpha (n)
= +\infty$ as well.  Therefore, no member of $\CH$ dominates the
function $f^E_\beta \in \CF$ and we conclude that \gd($\CF$)$\geq \kappa$
and thus \gd($\CF$)$=\kappa$. 

\noi {\bf \ref{U-class}.2:} We build an $\CF$ in the $\CU$-class such that \gb($\CF$)=2
and \gd($\CF$)$=\kappa$ for any  $\chi(\CU) \leq \kappa \leq $\gc.
	   
\noi For $X \in \CU$ such that $X \sub E$ and any $\alpha < \kappa$, let

\[ f^X_\alpha(n)= h(next(n,[A^c_\alpha \cap P(X)] \cup X)) + |X^c \cap n| \]

\noi and put $\CF = \{f^X_\alpha : X \sub E, X \in \CU, \alpha < \kappa
\}$.  For any $f^X_\alpha \in \CF$ and $x \in X$,we have $f^X_\alpha (x)
\leq h(x) + x$ and therefore $g$ and $h'(n)=h(n)+n$ witness that $\CF$
belongs to the $\CU$-class. 

We first show that the bounding number is 2.  No $f^X_\alpha$ itself is
unbounded in $\CF$ since choosing $Y \in \CU$ such that $X \setminus Y$
is infinite, we get for each $N$ and $n$ large enough

$\begin{array}{ll}
 f^X_\alpha(n) & =h(next(n,[A^c_\alpha \cap P(X)] \cup X)) + |X^c \cap n| \\
               &  \leq h(next(n,[A^c_\alpha \cap P(Y)] \cup Y)) + |X^c \cap n| \\
               &  \leq h(next(n,[A^c_\alpha \cap P(Y)] \cup Y)) + |Y^c \cap n| - N\\
               &= f^Y_\alpha (n) -N
\end{array}$

\noi and therefore $f^X_\alpha \prec f^Y_\alpha$. 

However, we claim that for any $\alpha \neq \beta$, the pair
$\{f^E_\alpha, f^E_\beta \}$ is unbounded in $\CF$.  To verify this, we
consider any $f^X_\gamma \in \CF$, without loss of generality $\alpha
\neq \gamma$.  For any $x \in A^c_\gamma \cap P(X) \setminus A^c_\alpha$
which is infinite, we have $f^X_\gamma (x) \leq h(x) + x$ and
$f^E_\alpha (x) \geq h(x+1)$.  But as $\limsup_n h(n+1) - [h(n) + n] = +
\infty$ by assumption we get $\limsup_n f^E_\alpha (n) - f^X_\gamma (n)
= +\infty $ as well.  The fact that \gd($\CF$)$=\kappa$ is proved as in
the previous example. 

\noi {\bf \ref{U-class}.3:} We build an $\CF$ in the $\CU$-class such that
 \gb($\CF$)$=\omega$ and \gd($\CF$)$=\kappa$ where  $\chi(\CU) \leq \kappa \leq $\gc.

\noi For any $a \in [\kappa]^{<\omega}$ and $X \in \CU$ such that $X \sub E$ we let

\[ f^X_\alpha(n)= h(next(n,[\bigcap_{\alpha \in a} A^c_\alpha \cap P(X)]
                     \cup X)) + |X^c \cap n| \]
 
\noi and put $\CF = \{f^X_\a : X \sub E, X \in \CU, a \in
[\kappa]^{<\omega} \}$.  Observe that for any such $f^X_a$ and $x \in X
\sub E$, $f^X_a (x) = h(x) + |X^c \cap x| \leq h(x) + x $
 and therefore $g$ and $h'(n)=h(n)+n$ again witness that $\CF$ belongs
to the $\CU$-class. 

Our first task is to show that $\CF$ is directed and therefore
\gb($\CF$)$\geq \omega$.  But given $f^X_a$ and $f^Y_b$, put $c=a \cup
b$ and choose $Z \in \CU$ such that $Z \sub X \cap Y$ and $X \cap Y
\setminus Z$ is infinite; then $f^Z_c \succ f^X_a, f^Y_b$.  To show now
that \gb($\CF$)$\leq \omega$, choose $A \in [\kappa]^\omega$ and we
prove that the collection $\CH= \{f^E_\alpha : \alpha \in A\} \sub \CF$
is unbounded in $\CF$.  So let us fix $f^X_a \in \CF$ and choose $\beta \in A \sm
a$.  Now the set $[\bigcap _{\alpha \in a } A^c_\alpha \cap P(X)] \sm
A^c_\beta$ is infinite as otherwise we would obtain $P(X) \cap A_\beta
\sub^* \bigcup_{\alpha \in a} A_\alpha$, and as $P(X) \cap A_\beta$ is
infinite $A_\beta$ would have infinite intersection with some $A_\alpha$
contradicting that $\CA$ is an almost disjoint family.  But now for $x
\in [\bigcap _{\alpha \in a } A^c_\alpha \cap P(X)] \sm A^c_\beta$,

\[ f^X_a (x) = h(next(x,[\bigcap_{\alpha \in a} A^c_\alpha \cap P(X)] \cup
     X))+ |X^c \cap x| \leq h(x)+ x, \]
\[  f^E_\beta (x) = h(next(x, A^c_\beta \cap P(E)] \cup E))
     + |E^c \cap x|  \geq h(x+1).\]

\noi As $\lim_n h(n+1)-[h(n)+n] = + \infty$, we get that $\limsup_n
f^E_\beta (n) - f^X_a (n) = + \infty$ as well and $\CH$ is indeed
unbounded in $\CF$. 

The verification the the dominating number is $\kappa$ is again very
similar to the first example. 

This completes the proof of Proposition \ref{U-class}. \qed

\begin{cor}
For any  $\lambda \in \{1,2,\omega\}$
and \gu \/ $ \leq \kappa \leq $\gc, there is a family $\CF$ in the
$\CU$-class such that \gb($\CF$)$=\lambda$ and \gd($\CF$)$=\kappa$. 
\end{cor}

\medskip

The next problem is whether we can construct a family $\CF$ in the
$\CU$-class with an uncountable bounding number.  We show
that this requires a $P$-point and therefore, in view of Shelah's
consistency result \cite{Shelah1} that there might be no such $P$-points,
we cannot construct such families in ZFC alone. 

\begin{prop}\label{P-point}
 If $\CF$ in the $\CU$-class has an uncountable
bounding number, then there is a finite-to-one map $m$ such that
$m(\CU)$ is a P-point.
  \end{prop}

\proof Fix functions $g$ and $h$ witnessing that $\CF$ belongs to the
$\CU$-class and define a sequence of integers such that $\pi_0=0$ and
more generally such that $g(\pi_{n+1}) > h(\pi_n)$.  We may assume
without loss of generality that $E= \bigcup_n [\pi_{2n}, \pi_{2n+1}) \in
\CU$ as the other case is analogous.  Now for any $X \in \CU$, if $X
\sub E$, any $f \in \CF$ with $g(next(-,X)) \leq^* f$ must satisfy

\[S(f)= \{n : f(n) \leq h(n) \} \sub^* T(X)=
    \bigcup\{[\pi_{2n-1},\pi_{2n+1}): X \cap [\pi_{2n},\pi_{2n+1}) \neq
                   \emptyset \}          \]

Now define a map $m \in \oo$ by $m"[\pi_{2n-1},\pi_{2n+1})=n$ and
consider $\CV=m(\CU)$.  Certainly $\CV$ is a (non principal) ultrafilter
as $m$ is finite-to-one.  To show it is actually a $P$-point, let $\{Y_n :
n \in \omega\} \sub \CV$ be given and consider the sets $X_n=
m^{-1}(Y_n) \cap E \in \CU$.  Since we are assuming that the bounding
number of $\CF$ is uncountable, fix a function $f \in \CF$ such that
$g(next(-,X_n)) \leq^* f$ for each $n$.  Therefore $S(f) \sub^* T(X_n)$
for each $n$ and thus $m(S(f)) \sub^* m(T(X_n)) \sub Y_n$.  Since
moreover $S(f) \in \CU$ and therefore $m(S(f)) \in m(\CV)$, the proof is
complete.  \qed

Under the existence of $P$-points or more generally $P_\kappa$ points,
one can easily construct members of the $\CU$-class with bounding number
$\kappa$ by fixing some $g \in \oo$ and defining $\CF=\{g(next(-,X)): X
\in \CU \}$.  Thus in general we have:

\begin{prop}\label{P-kappa}
There is a P$_\kappa$ ultrafilter if and only if there if a family $\CF$ in the
 $\CU$-class with bounding and dominating number $\kappa$.
\end{prop}

We can also deduce from the proof of Proposition \ref{P-point} that
\gd($\CF$)$\geq$\gu \/ for any family $\CF$ in the ultrafilter class, but I
do not know if \gd($\CF$)$\geq \chi(\CU)$ whenever $\CF$ belongs to the
$\CU$-class with witness $\CU$.

\section{Bounded Families}

Let $\CF$ be a bounded family and let

\[ \CF \!\! \downarrow = \{g \in \oo: (\fa f \in \CF)\;  f \prec g \} \]

\noi Certainly $\dCF$ is nonempty as $\CF$ is bounded and the pair $(\CF,\dCF)$ forms
a {\em gap\/} in the sense that there is no $h \in \oo$ such that
 \[(\fa f \in \CF)(\fa g \in \dCF)\; f \prec h \prec g .\]

To make a first distinction between bounded families, we make the
following definition. 

\begin{defn}

\[ \begin{array}{ll}
\mbox{\dgb}(\CF) &= \min \{ | \CH | : \CH \sub \dCF
                        \mbox{ and $\CH$ is unbounded in } \dCF \\
                 & \mbox{ in the reverse order } \} \\
                 &= \min\{\mid \CH \mid : \neg  (\ex g \in \dCF)(\fa h \in \CH)
                            \;  g \prec h \} 
\end{array} \]

\end{defn}

\noi We loosely call \dgb$(\CF)$ \/ and for that matter $\dCF$ depending
on the context the upper bound of $\CF$ and we will classify the
families according to this cardinal \dgb$(\CF)$ which takes either the
value 1 or an infinite regular cardinal; notice that the value 2 cannot
occur here.  Observe also that if $\CH \sub \dCF$ is unbounded in $\dCF$
in the reverse order as above, then the pair $(\CF,\CH)$ is also a gap. 
Much work has been done on gaps $(\CF,\CH)$ for which both $\CF$ and
$\CH$ are linearly ordered by $\prec$; in particular gaps $(\CF, \CG)$
for which \gb($\CF$)=\gd($\CF$).  Such gaps are usually qualified as
$($\gb$(\CF)$, \dgb$(\CF))$ gaps.  Here we will work in a more general
situation.

\subsection{Bounded families with a countable upper bound}

Unbounded families have much influence on the bounded ones; we can use
the results of \S3 to construct families with countable upper bounds and
various bounding and dominating numbers. 

\begin{prop} \label{unb}
There are families $\CF$ with countable upper bounds, that is \dgb$(\CF)=1$ or 
$\omega$, such that:
\begin{enumerate}
\item \gb($\CF$)=\gb \/ and \gd($\CF$)=\gd.
\item  \gb($\CF$)=1 or 2 and \gd($\CF$)=\gc.
\item  \gb($\CF$)=1, 2 or $\omega$ and \gu $\leq$ \gd($\CF$) $\leq$ \gc.
\end{enumerate}
\end{prop}

\proof The goal of the proof is to build families $\CF$ with the same bounding
and dominating number as the families from \S3; 
we fix for our constructions the functions
$g(n)=n^2$ and more generally for $\ell \in \omega$

\[ g_\ell(n)=n^2-\ell \log n  \mbox{ or } g_\ell(n)=n^2 - \ell
           \mbox{ depending on the context.} \]

\noi We shall build gaps $(\CF, \{g_\ell: \ell \in \omega \})$ giving us
familes $\CF$ with upper bounds 1 or $\omega$, depending on which
collection $\{g_\ell: \ell \in \omega\}$ one chooses, and with the
apropriate bounding and dominating numbers. 

Observe first that irrespective of the collection we choose, we have 
 \[ g_{\ell +1}(n+1) \geq g_\ell(n) \]
 for each $n$ and $\ell$; this will make our verifications easier. 
Now if $\CH$ is any unbounded family and $h \in \CH$, let

\[ f_h(n) = g_m(n) \mbox{ if } h(m-1) < n \leq h(m)  \]

\noi and put $\CF(\CH)= \{f_h: h \in \CH \}$. In this case we have :

\begin{claim}
\gb$(\CH)$=\gb$(\CF(\CH))$  and \gd$(\CH)$=\gd$(\CF(\CH))$. 
\end{claim}

\proof It suffices to prove that for all $h_1,h_2 \in \CH$, we have

\[h_1 \leq^* h_2 \mbox{  iff  } f_{h_1} \leq^* f_{h_2} .\]

\noi To verify this, suppose first that $h_1(m) \leq h_2(m)$ for all
$m \geq M$ and fix $n \geq h_2(M)$.

\noi Choose first $m$ such that
\[  h_1(m-1) < n \leq h_1(m)   \]
\noi and $\ell$ such that
\[h_2(\ell -1) < n \leq h_2(\ell). \]

\noi Observe that we must have $\ell \leq M$ and thus

\[ f_{h_1}(n) = g_m(n) \leq g_\ell(n) = f_{h_2}(n). \]

Suppose now for the other direction that $f_{h_1}(n) \leq f_{h_2}(n)$
for all $n \geq h_2(N)$ and fix $n \geq N$; we show that $h_1(n) \leq h_2(n)$.
But if for the sake of a contradiction we have $h_2(n) < h_1(n)$, pick
$\ell \geq n+1$ such that $h_2(\ell -1) < h_1(n) \leq h_2(\ell).$
Then
\[ f_{h_1}(h_1(n))=g_n(h_1(n)) \]
\noi and
\[ f_{h_2}(h_1(n))=g_\ell(h_1(n)) \leq g_{n+1}(h_1(n)) < g_n(h_1(n)) \]
and we obtain the desired contradiction. This proves the claim.

The Proposition is now proved by replacing $\CH$ by the apropriate
families of \S3. Actually, to obtain $\CF(\CH) \subseteq \oo$, we should 
first replace the families $\CH$ by $\CH'=\{h(n)+n: h \in \CH \}$ for
example to ensure that we have strictly increasing functions; observe that 
this does not affect the bounding and dominating number.  \qed

\medskip

There is however more than just reflecting unbounded families to bounded
ones, indeed let us see how close we are.  Let $\CF$ be a family of
functions and $\{g_n : n \in \omega \}$ a collection such that $ \CF
\prec g_{n+1} \leq^* g_n$ for each $n$, and assume without loss of generality that
$g_{n+1}(k) + 1 \leq g_n(k)$ for each $k$ and $n$.  For $f \prec \{g_n :
n \in \omega\}$,
 we define
 \[ h_f(n) = \max \{k : g_n(k) \leq f(k) \} \]

\noi and put $\CH(\CF)=\{h_f : f \in \CF\}$. The following proposition, due
to Rothberger, shows that unbounded families are always involved somehow.

\begin{prop}(Rothberger) \label{Roth}
  The pair $(\CF,\{g_n:n \in \omega\})$ is  a gap if and only if
$\CH(\CF)$ is an unbounded family.
\end{prop}

\proof   Suppose first that the family $\CH(\CF)$ is bounded, say by $h$; 
we might as well assume that $n<h(n)<h(n+1)$ for each $n$. Define a function $p$
by:
\[ p(j)=g_m(j) \mbox{ where $m$ is the smallest integer such that $h(m+1)>j$.} \]

\noi As $j$ increases, $m$ increases as well and therefore $p \prec g_m$ for each 
$m$. 
\noi Now for any $f \in \CF$, and therefore for $h_f \in \CH$, choose $N$ large
enough so that
\[ (\fa n \geq N)\;  h_f(n) < h(n). \]

\noi Hence for all $m \geq h(N)$, if we let $\ell \geq N$ be as large as possible
such that $m \geq h(\ell)$, we obtain:
\[ m \geq h(\ell) > h_f(\ell) \]

\noi and therefore

\[f(m) < g_\ell(m) =p(m). \]

\noi We conclude that $\CF \prec p \prec \{g_n:n \in \omega \}$ and thus
the pair $(\CF, \{g_n:n \in \omega\})$ is not a gap.

For the other direction, since we have the implication
\[f \leq^*f' \implies h_f \leq^* h_{f'},\]
 we conclude readily that
$\CH(\CF)$ is bounded if the pair  $(\CF, \{g_n:n \in \omega\})$ is not
a gap.   \qed

\begin{cor} \label{d}
\gd$(\CF)\geq$ \gb \/ for any $\CF$ with countable upper bound.
\end{cor}

Since $f \leq^*f' \implies h_f \leq^* h_{f'}$, we also obtain

\begin{cor}
\gd$(\CH(\CF)) \leq$ \gd$(\CF)$ and \gb$(\CF) \leq $ \gb$(\CH(\CF))$ unless
\gb$(\CH(\CF))=1$ in which case \gb$(\CF) \leq 2$.
\end{cor}

This allows us to extend Proposition \ref{unb} as follows.

\begin{prop}\label{gen} Let $\CH$ be any unbounded family and $\lambda
\leq$\gb \/ a regular (infinite) cardinal. Then there is a family $\CF$
with countable upper bound such that

\[ \mbox{\gb}(\CF)=\min \{\lambda,\mbox{\gb}(\CH)\} \mbox{and }
               \mbox{\gd}(\CF)=\mbox{\gd}(\CH).    \]
\end{prop}

\proof   To simplify the calculations, we fix the functions
$g_k(n)=n^n-kn$ for $k \in \omega$ and an increasing sequence of sets
$\la X_\alpha: \alpha < \lambda \ra$ such that $X_\beta \sm X_\alpha$ is 
infinite whenever $\alpha < \beta$; this is guaranteed by $\lambda \leq$ \gb.

\noi Without loss of generality, we may assume that each $h \in \CH$ is 
strictly increasing, that $h(n)>n$ for each $n$ and that the range is included
in $X_0$. Now for $h \in \CH$ and $\alpha < \lambda$, define

\[f_{h,\alpha}(n)=g_m(last(n,X_\alpha))+|X_\alpha \cap n|
             \mbox{ where } h(m-1) < n \leq h(m) \]

\noi and put $\CF=\{f_{h,\alpha} : h \in \CH, \alpha < \lambda \}$.  As
$\CH(\CF)=\CH$, we conclude from Corollary 4.6 that $( \CF, \{g_k:k \in
\omega \} )$ is a gap and that \gb$(\CF) \leq$\gb$(\CH)+1$ and
\gd$(\CF)\geq$\gd$(\CH)$. 

\begin{claim}
\gb$(\CF)\leq \lambda$.
\end{claim}

\proof Fix $h \in \CH$ and let $\CS=\{f_{h,\alpha}: \alpha < \lambda\}$. 
We show that $\CS$ ($\subseteq \CF$) is unbounded in $\CF$. Indeed, fix any
$h' \in \CH$ and any $\alpha < \lambda$ and consider any $\beta$,
$\alpha < \beta < \lambda$; we claim that $f_{h,\beta}(n) \geq f_{h',\alpha}(n)$ 
for infinitely many $n$, indeed on almost all $x \in X_\beta \sm X_\alpha$. For fix
such an $x$,  if $h'(m-1) < x \leq h'(m)$, 
\noi then 
\[ \begin{array}{ll}
  f_{h',\alpha}(x) & = g_m(last(x,X_\alpha)) + |X_\alpha \cap x| \\
                   & \leq g_m(x-1) + x \\
                   & = (x-1)^{x-1} -m(x-1) + x
\end{array} \]

\noi and if $h(\ell-1)< x \leq h(\ell)$  then
\[ \begin{array}{ll}
  f_{h,\beta}(x)  & = g_\ell(last(x,X_\beta)) + |X_\beta \cap x| \\
                   & \geq g_\ell(x) = x^x - \ell x.
\end{array} \]

\noi As $\ell,m \leq x$, we get $f_{h,\beta}(x) \geq f_{h',\alpha}(x)$ for almost
all such $x$'s. This proves the claim.

\begin{claim}
\gb$(\CF)\geq \min \{\lambda, \mbox{\gb}(\CH)\}$.
\end{claim}

\proof Let $\CS \sub \CF$, $|\CS| < \min\{\lambda, \mbox{\gb}(\CH)\}$, and fix 
$\zeta < \lambda$, $\CT \sub \CH$ such that
\[\CS \sub \{f_{h,\alpha}: h \in \CT, \alpha < \zeta \}
            \mbox{ and } |\CT| < \mbox{ \gb}(\CH). \]

\noi Therefore choose an $h' \in \CH$ such that $h <^* h'$ for any $h \in \CT$ and
we show that $f_{h,\alpha} <^* f_{h',\zeta}$ for all $h \in \CT$ and $\alpha 
< \zeta$, and thus $\CS$ is bounded in $CF$.

Choose first $N \in X_\alpha$ such that $X_\alpha \sm N \sub X_\zeta$ and fix 
$n \geq N$; if $m$ is such that
\[h'(m-1) < n \leq h'(m), \]
\noi and $\ell$ such that
 \[ h(\ell -1) < n \leq k(\ell ) \]
\noi we obtain, with $x=last(n,X_\zeta)$,

$\begin{array}{ll}
 f_{h,\alpha}(n) & = g_\ell(last(n,X_\alpha)) + |X_\alpha \cap n| \\
                 & \leq g_\ell(last(n,X_\zeta)) + |X_\alpha \cap n| \\
                 & = g_\ell(x) + |X_\alpha \cap n | \\
                 & = x^x - \ell x + |X_\alpha \cap n|
\end{array}$

\noi and

$\begin{array}{ll}
f_{h',\zeta}(n) & = g_m(last(n,X_\zeta)) + |X_\zeta \cap n| \\
                & = g_m(x) + |X_\zeta \cap n| \\
                & = x^x - mx + |X_\zeta \cap n|.
\end{array}$

\noi But $m \leq \ell$ (for $n$ large enough) and as $X_\zeta \sm
X_\alpha$ is infinite, we get $f_{h,\alpha}(n) < f_{h',\zeta}(n)$ for
almost all $n$, in fact $f_{h,\alpha} \preceq f_{h',\zeta}$.  This
proves the claim. 

\medskip

Finally, as we already know that \gd$(\CF) \geq$ \gd$(\CH)$, we must
show the reverse inequality.  But $\CF$ is generated by $\lambda \times
$\gd$(\CH)=$\gd$(\CH)$ functions, and the proof is complete. \qed

This gives an idea of what can be done in terms of bounded families with
countable upper bound, they all involve unbounded
families by proposition \ref{Roth}, but this is only very partial
information and a lot of freedom remains.

\subsection{Families with upper bound $\omega_1$}

One of the surprising construction in ZFC is a gap build by Hausdorff
which has bounding and dominating number $\omega_1$.  Lusin build one
with bounding number 1 and dominating number $\omega_1$; it is this
construction that we will adapt to produce gaps with various bounding
and dominating numbers.  Although in both Hausdorff's and Lusin's
construction the upper bound is at most $\omega_1$, I do not know if
could be $\omega$; if \gb$> \omega_1$, they certainly cannot. 

\begin{prop} \label{omega1}
For each $\lambda \in \{1,2,\omega\}$ and $\omega_1 \leq
\kappa \leq$ \gc, there is a family $\CF$ with upper bound  at most $\omega_1$
such that \gb($\CF$)=$\lambda$ and \gd($\CF)=\kappa$.
\end{prop}

\proof

We first build $\{f_\alpha : \alpha < \omega_1 \}$ and $\{g_\alpha : \alpha < \omega_1 \}$
 such that:

\begin{description}

\item[1:] $(\fa \alpha < \beta)\;  f_\alpha + id \prec g_\beta \prec g_\alpha$,
         where $id$ is the identity function $id(n)=n$.

\item[2:] $(\fa \alpha)(\fa a \in [\omega_1 \sm \{\alpha\}]^{< \omega}) \; 
   \limsup_n f_\alpha(n) - \max_{\gamma \in a} \{f_\alpha(n) + n\} = + \infty .$

\item[3:] $(\fa \alpha)(\fa n)\;  f_\alpha(n) \leq g_\alpha(n)$.

\item[4:] $(\fa \alpha < \beta )(\ex n) \;  f_\alpha(n) > g_\beta(n)$.

\end{description}

Let us first observe that this construction, essentially due to Lusin, 
will give us a gap.

\begin{claim} 
The collection ${\la } \{f_\alpha : \alpha < \omega_1 \},
\{g_\alpha : \alpha < \omega_1 \} {\ra }$ is a gap. 
\end{claim}

\proof  Suppose on the contrary that
$ \{f_\alpha : \alpha < \omega_1 \} \prec h \prec \{g_\alpha : \alpha < \omega_1 \}$
 for some function $h$. Choose $X \in [\omega_1]^{\omega_1}$ and $n$ so that

\begin{description}

\item[a:] $(\fa \alpha, \beta \in X) \;  f_\alpha \rest n = f_\beta \rest n
          \mbox{ and } g_\alpha \rest n = g_\beta \rest n$.

\item[b:] $(\fa m \geq n) \;  f_\alpha(m)  \leq g_\beta (m)$.

\end{description}

\noi Thus 
 $(\fa \alpha < \beta \in X)(\fa k)$
         \[ \begin{array}{lll}
             f_\alpha(k) & =f_\beta(k) \leq g_\beta(k) &  \mbox{ if } k < n  \\
             f_\alpha(k) & \leq g_\beta (k)   &  \mbox { if } k \geq n
           \end{array} \]

\noi But this contradicts requirement 4. \qed

\vspace{.1in}

If we can accomplish this construction, we put $\CF_1 = \{f_\alpha :
\alpha < \omega_1\}$, $\CG = \{g_\alpha : \alpha < \omega_1 \}$ and we
get a gap $(\CF,\CG)$ with \gb($\CF_1)=1$ and \gd($\CF_1)=\omega_1$.  
Choosing functions $0 \prec h_n \prec h_{n+1} \prec id$ and using
$\CF_2 = \{f_\alpha + h_n:  \alpha < \omega_1, n \in \omega \}$,
 we obtain a family with \gb($\CF_2)=2$ and \gd($\CF_2)=\omega_1$. 
Finally, we let
 $\CF_\omega = \{\max_{\alpha \in a}\{f_\alpha \}+ h_n : a \in [\omega_1]^{< \omega}, 
n \in \omega \}$ 
we obtain a family with \gb($\CF)=\omega$ and
\gd($\CF)=\omega_1$.
 To obtain familes with various dominating number, fix for example
$f_0$ and choose a set $X=\{x_n:n \in \omega \}$ such that $f_0(x_{n+1}) > f_0(x_n)
+x_n$ and let $\CA=\{A_\alpha: \alpha < \kappa \} \sub \CP (X)$ an almost disjoint family.
Assume further that we actually have $0 \prec 2h_n \prec 2h_{n+1} \prec id$.
Now for $\beta < \kappa$, define

\[f^\beta_0(n)=  \max\{f_0(n),f_0(last(n,A_\beta))+ \frac{1}{2}last(n,A_\beta)\} .\]

\noi Notice that for $\beta \neq \beta'$, if $x_{n+1} \in A_\beta \sm a_{\beta'}$, 
\[\begin{array}{ll}
  f^\beta_0(x_{n+1}) & = f_0(x_{n+1}) + \frac{1}{2}x_{n+1} \\
  f^{\beta'}_0(x_{n+1}) & \leq \max\{f_0(x_{n+1}),f_0(x_n) + x_n\} \\
                        & = f_0(x_{n+1})
\end{array} \]

\noi and therefore $f^\beta_0(x_{n+1}) - f^{\beta'}_0(x_{n+1}) \geq x_{n+1}$ and hence
for each $m$ 
\[ \limsup_k f^\beta_0(k) - [f^{\beta'}_0(k)+h_m(k)] = + \infty . \]

\noi We can then replace $f_0$ in the above families by $\{f^\beta_0: \beta < \kappa\}$
to obtain families with dominating number $\kappa$.

\vspace{.2in}

The construction proceeds by induction on $\alpha$, that is we start with

\[ f_0(n)=n \mbox{ and } g_0(n)=n^2 \]

\noi Now assume that we have already defined the functions 
$\{f_\xi: \xi \in \alpha \}$ and $\{g_\xi: \xi \in \alpha \}$
such that:

\begin{description}

\item[2.1:] $(\fa \beta < \gamma < \alpha)\;  f_\beta + id \prec g_\gamma \prec g_\beta$.

\item[2.2:] $(\fa \beta < \alpha)(\fa a \in [\alpha \sm \{\beta\}]^{< \omega})
           \;   \limsup_n f_\beta(n) - \max_{\gamma \in a} \{f_\gamma(n) + n \} = + \infty$.

\item[2.3:] $(\fa \beta < \alpha)(\fa n)\;  f_\beta(n) \leq g_\beta(n)$.

\item[2.4:] $(\fa \beta < \gamma < \alpha)(\ex n)\;  f_\beta (n) > g_\gamma(n)$.

\end{description}

\noi and we proceed to build $f_\alpha $ and $g_\alpha$ in countably many steps.
As $\alpha$ is countable, we list $\alpha \times [\alpha]^{< \omega}$
as $\{{\la  }\alpha_k, a_k {\ra  }: k \in \omega \}$, $\{f_\beta : \beta < \alpha \}$ as
$\{f^k : k \in \omega \}$ and $\{g_\beta : \beta < \alpha \}$ as
$\{g^k : k \in \omega \}$.

\noi At stage $N$, suppose that we have $f_\alpha \rest n$ and $g_\alpha \rest n$
for some $n$,  such that:

\begin{description}

\item[3.1:] $(\fa k < N)(\ex m < n)\;
  f_{\alpha_k}(m) - \max_{\gamma \in (a_k \cup \{\alpha\})\sm \{\alpha_k\}}
                         \{f_\gamma(m) + m \} \geq N$,

\item[3.2:] $(\fa k < N)(\ex m < n) \;
              f_\alpha(m) - \max_{\gamma \in a_k} \{ f_\gamma(m) + m \} \geq N$,

\item[3.3:] $(\fa m < n) \; f_\alpha(m)  \leq g_\alpha(m)$,

\item[3.4:] $(\fa k < N)(\ex m<n) \; f^k(m) > g_\alpha(m)$.

\end{description}

\noi We will also ensure that for $m \geq n$

\noi {\bf 3.5:} $\max_{k<N} \{f^k(m)+m\} + 2N \leq g_\alpha(m) + N \leq
\min_{k<N} \{g^k(m)\}.$

\noi This will help satisfy 2.1.  Requirements 3.1 and 3.2 will ensure
2.2, 3.3 will give 2.3 and 3.4 gives 2.4. 

For the construction at stage $N$, first choose $m_0 > n$ such that:

\begin{enumerate}

\item $f^N(m_0) - \max_{k<N} \{f^k(m_0)+m_0 \} > N$.

\item $f^N(m_0) > g_\alpha(n-1)$.

\end{enumerate}

\noi Then we define, for $n \leq m \leq m_0$,
\[ \begin{array}{ll}
    g_\alpha(m) & = \max_{k<N} \{g_\alpha(n-1), f^k(m) + m + N\} \\
    f_\alpha(m) &= f_\alpha(n-1)
   \end{array} \]

\noi This fulfills 3.4 as $g_\alpha(m_0) < f^N(m_0)$.

Now given $m_i$ for $i<N+1$, choose $p_i > m_i$ such that:

\begin{enumerate}

\item $f_{\alpha_i}(p_i) - \max_{\gamma \in a_i \sm \{\alpha_i\}}
                  \{f_\gamma(p_i) + p_i \} \geq N+1$.
\item $\max_{\gamma \in a_i \sm \{\alpha_i\}} \{f_\gamma (p_i) \} \geq f_\alpha (m_i)$.

\end{enumerate}

\noi Then we define, for $m_i \leq m \leq p_i$,
\[ \begin{array}{ll}
   f_\alpha(m) &= f_\alpha(m_i) \\
   g_\alpha(m) &= \max_{k<N}\{g_\alpha(m_i),f_\alpha(m)+N,f^k(m)+m+N\}.
\end{array} \]

\noi This handles 3.1.

Finally choose $n' > p_N$ large enough so that for all $m \geq n'$, we have

\[\begin{array}{ll}
  g_\alpha(p_N) & \leq \max \{ \max_{\gamma \in \cup_{k<N+1} a_k} \{f_\gamma(m)+m\},
                               \max_{k<N+1} \{ f^k(m)+m \}  \}  \\
                & \leq \min_{k<N} \{g^k(m)\}  - 3(N+1) 
   \end{array} \]

\noi and define

\[ \begin{array}{ll}
   f_\alpha(m) &= f_\alpha(p_N)    \mbox{ for } p_N < m < n' \\
               &= \max_{\gamma \in \cup_{k<N+1} a_k} \{f_\gamma(n')+n'\} + N + 1 \\     
               &  \hspace{.6in} \mbox{ for } m=n' \\
   g_\alpha(m) &= \max_{k<N} \{ f_\alpha(m) + N, f^k(m) + N \} \mbox{ for } p_N < m \leq n'
   \end{array} \]

\noi This satisfies 3.2 and observe that we are able to keep our promise
3.5.  This completes the construction and proves the Theorem \ref{omega1}.  \qed

As far as uncountable bounding number is concerned, a Hausdorff gap provides a family
$\CF$ with upper bound at most  $\omega_1$, \gb$(\CF)$=\gd$(\CF)=\omega_1$. I do not know
if there is always such a family $\CF$ with large dominating number, say \gd$(\CF)=$\gc \/
for example.

\subsection{Families with upper bound \gb}

In view of Rothberger's result and the fact that the smallest size of an unbounded family
in $\oo$is \gb, it is not at all surprising that this cardinal has a role to play in
bounded families. We have the following result.

\begin{prop}\label{b}
For any $\lambda \in \{1,2,\omega\}$, and  $\lambda \leq \kappa \leq$ \gc, there is a 
family $\CF$ with upper bound \gb \/  such that \gb($\CF$)$=\lambda$ and
\gd($\CF$)$=\kappa$.
\end{prop}

\proof We provide a general construction which will work for all values of $\lambda$
and $\kappa$.

Fix an increasing unbounded family ${\la }h_\alpha: \alpha <
\mbox{\gb}{\ra }$ and let $f(n)=n^2$.  Now for $\ell \in \omega$ define
\[ f_\ell(n)= n^2 + \ell \log (n), \mbox{ and thus } f_\ell \prec
f_{\ell +1} \]

\noi and for $\alpha <$ \gb \/ put

\[g_\alpha(n)=f_m(n)=n^2 + m \log (n) \mbox{ if } h_\alpha(m-1) < n \leq h_\alpha (m). \]

\noi These functions are technically not in $\oo$ because of the $\log$
function, but one could easily take instead the smallest
integer greater than or equal to these values. 
Notice that

\[f_\ell \prec f_{\ell +1} \prec g_\beta \prec g_\alpha \mbox{ for all } \alpha < \beta
                  \mbox{ and } \ell. \]

\begin{claim}
For all $X \in \so$, the pair ${\la  }\{f_\ell \rest X: \ell \in \omega \}, \{g_\alpha \rest X: 
\alpha < \mbox{\gb} \} {\ra  }$ is a gap.
\end{claim}

\noi {\bf Proof of the claim:} Suppose otherwise that there is a function $h$ such that

\[ (\fa \ell)(\fa \alpha)\;   f_\ell \rest X \prec  h \prec g_\alpha \rest X . \]

\noi Then we define

\[  p(n)= \min \{ x \in X: (\fa y \in X \sm  x) h(y)>f_n(y) \} \]

\noi It now suffices to show that $h_\alpha \leq^*p$ for each $\alpha$ to obtain
a contradiction.  

\noi But fix $N$ large enough so that
\[(\fa x \in X) \; x \geq N \implies h(x) < g_\alpha(x) \]
\noi So for each $n$ with $p(n) \geq N$ we must have $p(n) \geq h_\alpha(n)$ as well;
indeed, if $x=p(n)<h_\alpha(n)$, we get 

\[ h(x) = h(p(n)) > f_n(x). \]

\noi Now choose $m \leq n $ such that $h_\alpha(m-1) < x \leq h_\alpha (m)$,
then

\[g_\alpha(x) = f_m(x) \leq f_n(x) \]

\noi and therefore $g_\alpha(x) < h(x)$, a contradiction. This proves the claim. \qed

\vspace{.2in}

Now let $\CA = \{A_\alpha : \alpha < \kappa \}$ be an almost
disjoint family and for each $\alpha < \kappa$ and
 $\ell \in \omega$, let

\[ f'_{\alpha, \ell} (k) = f_\ell(last(n,A_\alpha))  \]

\noi Then certainly $f'_{\alpha, \ell} \prec f'_{\alpha, \ell+1}$. 
Further, if $\alpha \neq \beta$ and $\ell, k$ are given,  pick 
$n \in A_\beta \sm A_\alpha$, and thus

\[f'_{\beta, \ell}(n) = f_\ell(n) = n^2 + \ell \log (n) \]

\noi  but as  $last(n,A_\alpha) \leq n-1$ we obtain

\[ f'_{\alpha, k}(n) \leq (n-1)^2 + k \log (n-1) = n^2 -2n +1 + k \log (n-1) \]

\noi and hence $\limsup_n f'_{\beta, \ell}(n) - f'_{\alpha, k}(n) = + \infty$.

If we now let $\CF = \{f'_{\alpha, \ell} : \alpha < \kappa, \ell \in \omega \}$, we obtain 
a family with upper bound \gb, bounding number 2 and dominating  number $\kappa \times \omega
= \kappa$. 

On the other hand if we let $\CF = \{ \max\{f'_{\alpha, \ell} : \alpha \in a \} : a \in 
[\kappa]^{< \omega}, \ell \in \omega \}$, we obtain a family with again upper bound \gb
\/ but bounding number $\omega$ and dominating number $\kappa$.

Moreover, if we had used the functions $f_\ell (n) = n^2 +\ell$ instead,
then the family $\CF = \{f'_{\alpha, 0} : \alpha < \kappa \}$ would
constitute a family with upper bound \gb,\/ bounding number 1 and
dominating number $\kappa$. 

This completes the proof of Proposition \ref{b}. \qed

\vspace{.1in}

The obvious question now is whether we can have a family with upper bound \gb \/
and uncountable bounding number; we shall see that there is no such family in the 
Mathias model and hence such familes cannot be constructed in ZFC alone.

\section{Models with few families of functions}

We shall be interested in two forcing notions.

\begin{defn}~
 \bd

\item[Mathias forcing] ${\Bbb M}_1 =\{ \la   a,A\ra  : a \in
[\omega]^{<\omega},A $ is an infinite subset of $\omega$ disjoint from
$a \}$ equipped with the ordering
  \[ \la   a,A\ra   \leq \la   b,B\ra   \mbox{ iff } A \sub B, b \sub a,
             \mbox{ and } a\sm b \sub B. \]

\item[Matet forcing] ${\Bbb M}_2 = \{ \la   a,A\ra   : a \in [\omega]^{< \omega}, A $ is
     an infinite set of pairwise disjoint finite subsets of $ \omega \}$ equipped with the 
     ordering
     \[ \la   a,A\ra   \leq \la   b,B\ra   \mbox{ iff } b \sub a, a \sm b \]
\noi and members of  $A$ are finite unions of elements of $B$.
 \ed
\end{defn}

\noi We use ${\CM}_1$ and $\CM_2$ to denote the models obtained from a
model of CH by an $\aleph_2$ iteration with countable support of the
(proper) partial orders ${\Bbb M}_1$ and ${\Bbb M}_2$ respectively. 

It is known from \cite{Blass2} that $\CM_2$ satisfies \gu$<$\gothg \/
and hence by \cite{Laf1, Laf2} that the only unbounded (downward closed)
families in this model are from the 3 classes described in \S3. 
Further, $\CM_2$ satisfies \gb=\gu=$\aleph_1$ and \gd=\gc=$\aleph_2$. 
On the other hand, the Mathias model $\CM_1$ satisfies
\gb=\gd=\gu=\gc=$\aleph_2$. 

\noi One can modify Baumgartner's result that Mathias forcing preserves
towers to gaps with uncountable upper bound and bounding number and
extend it to Matet forcing as well:

\begin{prop}
If $\CM \models ZFC$ and {\bf G} is either ${\Bbb M}_1$ or ${\Bbb M}_2$ generic, then
any $(\omega_1,\omega_1)$ gap from $\CM$ remains a gap in $\CM[${\bf G}$]$.
\end{prop}

\noi The iteration lemmas of Shelah \cite{Shelah2} give us:

\begin{prop}
In $\CM_1$ or $\CM_2$, the only linearly orderd $(\lambda, \kappa)$ gaps have either
$\lambda \leq \omega$, $\kappa \leq \omega$ or else $\lambda=\kappa=\omega_1$.
\end{prop}

Here are therefore the families we get in $\CM_1$.  From Propositions
\ref{unb} we get families $\CF$ with countable upper bound with
\gb$(\CF)=1,2$ or $\omega$ and \gd$(\CF)=\omega_2$.  By Proposition
\ref{gen}, we get families with countable upper bound, bounding number
$\omega_1$ or $\omega_2$ and dominating number $\aleph_2$ by using
$\CH=\oo$.  There are no such families with \gd$(\CF)=\omega_1$ by
Corollary \ref{d}.  For families with upper bound $\omega_1$ or
$\omega_2=$\gb, Propositions \ref{omega1}, \ref{b} and Hausdorff's
result provide a general context, although I do not know if $\CM_1$ has
a family with upper bound and bounding number $\omega_1$, and dominating
number $\aleph_2$; there is however an unbounded family $\CF$ in $\CM_1$
with \gb$(\CF)=\omega_1$ and \gd$(\CF)=\omega_2$.  The two Propositions
above 5.2 and 5.3 justify our remark of \S4.3 that no gaps with upper
bound \gb \/ has uncountable bounding number in $\CM_1$.  Indeed, a
standard argument would force such a family $\CF$ to reflect to some
$\CF_{\alpha}= \CF \cap \CM_1[${\bf G}$_\alpha]$ for some $ \alpha <
\omega_2$ where \gb$(\CF_\alpha)=$\dgb$(\CF_\alpha)=\omega_1$ and be
equivalent in this model to a linearly ordered $(\omega_1, \omega_1)$
gap.  Since this gap would be preserved to $\CM_1$, we obtain
\dgb$(\CF)=\omega_1$, a contradiction. 

In $\CM_2$ we have a little more:

\begin{prop}
In $\CM_2$ there are no $<^*$-increasing or $<^*$-decreasing chains  of
length $\omega_2$.
\end{prop}

\proof It suffices to prove that there are no increasing chains of size
$\omega_2$.  Let $\CF$ be such a chain.  If $\CF$ is unbounded, it would
have to belong to one of the 3 classes described in \S3.  Clearly $\CF$
cannot be dominating as \gb \/ would then have to be $\omega_2$ in this
model; $\CF$ cannot belong either to the $\CS$-class as
\gb$(\CF)=\omega_2 > 2$.  If finally $\CF$ would belong to the
$\CU$-class, then Proposition \ref{P-kappa} would provide us with a
$P_{\aleph_2}$-point which do not exist in $\CM_2$ by \cite{BlassLaf}. 
Therefore $\CF$ must be bounded and the above preservation results show
that \dgb$(\CF)$ is countable; then Rothberger's result, Proposition
3.2, produces an unbounded $<^*$-increasing chain $\CH(\CF)$ of size
$\omega_2$ which we have just showed does not exist.  \qed

This provides an alternative model to Theorem 3.1 of Shelah and Steprans
\cite{ShelahSteprans} showing the failure of Nyikos' axiom 6.5.  Indeed
the above shows that any family has an unbounded susbset of size at most
$\omega_1$ and since NCF holds in $\CM_2$ as it follows from
\gu$<$\gothg, we conclude that cof$(\omega^\omega/\CU)$ =\gd=$\aleph_2$
for all ultrafilters $\CU$ (see \cite{Blass0}). 

To summarize, we have the following bounded families in $\CM_2$.  For
$\lambda=1,2$ or $\omega$ and $\aleph_1 \leq \kappa \leq \aleph_2$,
Proposition \ref{unb} gives us $\CF$ with countable upper bound such
that \gb$(\CF)=\lambda$ and \gd$(\CF)=\kappa$.  For $\lambda =
\omega_1$, Proposition \ref{unb} again gives us $\CF$ countable upper
bound, bounding number $\omega_1$ and dominating number $\aleph_2$; as
\gu=$\aleph_1<$\gd=$\aleph_2$, we get a $P_{\aleph_1}$-point in $\CM_2$
and Proposition \ref{P-kappa} together with Proposition \ref{gen} give
us a family with countable upper bound, bounding number $\omega_1$ and
dominating number $\aleph_1$.  There are no such families with bounding
number $\omega_2$ as remarked above.  Now for families with upper bound
$\omega_1$, there those with bounding number $\lambda=1$,2 or $\omega$
and dominating number between $\lambda$ and $\omega_2$ by Propositions
\ref{omega1} and \ref{b} as \gb$=\omega_1$.  Hausdorff's result provides
one with bounding and dominating number $\omega_1$ and again I do not
know if there is one with (upper bound $\omega_1$) bounding number
$\omega_1$ and dominating number $\omega_2$.  There are no families with
upper bound $\omega_2$.

\end{document}